\documentclass{notices}
\usepackage{amsfonts,amssymb,amsmath,amscd}
\usepackage{url}

\title{
A source list to support DEI/EDI work in mathematical sciences
}

\author{
 Deborah Kent
  \affil{Reader in History of Mathematics, University of St Andrews, dk89@st-andrews.ac.uk 
  }
 \and
 Emilie Aebischer
  \affil{University of St Andrews, PhD candidate, Intellectual History, ea94@st-andrews.ac.uk
   }
\and
Stuart Neave
  \affil{University of St Andrews, PhD candidate, Modern History, san7@st-andrews.ac.uk
   }
  }

\begin{document}

\maketitle

Many mathematicians have recently become interested in considering matters of diversity, equity, and inclusion (DEI/EDI) in their professional work. For some, this is prompted by documents such as the \href{https://www.cde.ca.gov/ci/ma/cf/}{Mathematics Framework} proposed in 2021 by the California Department of Education, or the Quality Assurance Agency (QAA) 2022 revised subject benchmark for \href{https://www.qaa.ac.uk/the-quality-code/subject-benchmark-statements/subject-benchmark-statement-mathematics-statistics-and-operational-research }{Mathematics, Statistics, and Operational Research} (MSOR) in the United Kingdom, both of which stress the importance of considering DEI/EDI work for the discipline. The QAA benchmark states that ``Equality, diversity and inclusion (EDI) is essential for the health of MSOR, and it is important that the discipline encourages inclusivity and access to ensure learners are attracted from diverse backgrounds, that the curriculum and environment enable them to succeed in their studies, and that the subject is enriched by input from diverse practitioners.''

The source list presented in this article is available in full and linked here to \href{https://arxiv.org/archive/math}{arXiv}, is designed as a resource to support work in the education and practice of mathematics that involves DEI/EDI. The project developed alongside a UK network funded by the Issac Newton Institute to create an environment for discussing the place of DEI/EDI in mathematics curricula in higher education. 

The associated source list provides a variety of starting points for people interested in a range of concerns related to gender, race, ethnicity, (dis)ability, class, or sexuality that currently impact the study and practice of mathematical sciences at the post-secondary level. It makes no claims of completeness – many individual sources listed here could generate entire bibliographies of other relevant material – but is intended to provide a collection of references with a variety of perspectives that can serve as an introduction to relevant literature for mathematical practitioners. There is no explicit or implicit instruction to read these in any particular order, or to read every source on the list. Indeed, interested individuals will be in different institutional contexts with a varied range of interests and priorities related to DEI/EDI initiatives. This document is intended as a resource for anyone looking to deepen their own understanding of societal oppression, discrimination, inequity; how these things could impact the teaching and learning of mathematics; and ways instructors might redress the situation. It will necessarily involve individual efforts to implement resulting insights for a particular audience or approach. Each section of the reading list is introduced by a brief overview to assist readers with identifying materials of interest to them.
 
The selected sources are organised by scale, starting with general interest works that interrogate gender, race, ethnicity, (dis)ability, class, or sexuality in contemporary society. This is followed by articles and books that present research results about how matters of inequity, exclusion, and homogeneity surface in STEM educational contexts. Then come sources focused specifically on the past and present of interplay between these themes and the study and teaching of mathematics.  Each section includes a range of material organised alphabetically by author and not further categorised. Taken all together, these sources indicate some contours of current challenges that face mathematicians engaged in efforts to create diverse, equitable, and inclusive classroom environments. 

The second part of the list looks at visions for and approaches to creating better environments for all students of mathematics in higher education. It begins with a collection of sources that have a conceptual or theoretical orientation towards the challenges of creating more inclusive mathematics classrooms. The list then moves to consider approaches that specifically involve the history of mathematics as a tool to address the challenges faced by the community of mathematical practitioners. There are two primary threads of related research: one focuses on using primary source material to teach mathematical content, while the other draws on historical research to increase representation in mathematical sciences. 

\section*{General social/cultural context}

These general-interest books introduce a broad range of issues related to inequity and discrimination in contemporary society. One core theme is a call to recognize normative conditions that create discriminatory perceptions of “non-conforming” people.  Both Lennard J. Davis in the introduction to \textit{The Disability Studies Reader} and Fiona Kumari Campbell in a more technical tome, \textit{Contours of Ableism}, interrogate normativity in relation to perceptions of people with disabilities. Similar considerations of exclusionary norms appear in Reni Eddo-Lodge’s \textit{Why I’m No Longer Talking to White People About Race} and the introduction to Nikesh Shukla’s \textit{The Good Immigrant}. 
Another theme in this section is the recovery of neglected histories. David Olusoga’s \textit{Black and British: A forgotten history} counters the narrative that excludes Black people from British history before Windrush. In \textit{Racism, Class, and the Racialized Other}, Satnam Virdee presents a social history of Britain where race is a vital factor in the development of a working class. Selina Todd’s \textit{Snakes and Ladders: The Great British Social Mobility Myth} is another historical work that provides social context relevant to the matter of higher education in mathematics and other subjects. 
While most of these works aim at a general audience, there is variation in approach and tone. Some, like \textit{The Routledge Handbook of Epistemic Injustice} are considerably more technical than public-facing titles. On the more personal-development end of the spectrum, Layla Saad’s \textit{Me and white supremacy} encourages the reader to engage in self-reflective journalling alongside the text. Other works, like Ibram X. Kendi’s \textit{How to be an Antiracist} are more autobiographical. Similarly, Eileen Pollack’s \textit{The Only Woman in the Room} relies on her own experience as a mathematics student at Yale to illuminate past and present biases in the discipline.\\ 

\noindent Campbell F K (2009) Contours of Ableism: The Production of Disability and Abledness. Palgrave Macmillian, London\\

\noindent Collins P H, Bilge S (2016) Intersectionality. Polity, Cambridge\\
 
\noindent Cooper, B (2019) Eloquent Rage: A Black Feminist Discovers Her Superpower. St. Martin’s Press, New York\\

\noindent Davis, L J (2013) Introduction: Disability, Normality, and Power. In: Lennard J D (ed), The Disability studies Reader. Routledge, Abingdon, pp. 1-16\\

\noindent Eddo-Lodge R (2018) Why I’m no Longer Talking to White People about Race. Bloomsbury Publishing, London HERE\\

\noindent Kendi I X (2023) How To Be an Antiracist. The Bodley Head, London\\

\noindent Kidd I J, Medina J, Pohlhaus Jr G (2017) The Routledge Handbook of Epistemic Injustice. Routledge, London\\

\noindent Nash K (2022) Positively Purple: Build an Inclusive World Where People with Disabilities Can Flourish. Kogan Page, London\\

\noindent Olusoga D (2021) Black and British: A Forgotten History. Picador, London\\

\noindent Pérez C C (2019) Invisible Women: Data Bias in a World Designed for Men. Abrams Press\\

\noindent Pollack E (2015) The Only Woman in the Room: Why Science Is Still a Boys’ Club. Beacon Press, Boston\\

\noindent Saad L (2020) Me and White Supremacy: How to Recognize Your Privilege, Combat Racism, and Change the World. Quercus, London\\

\noindent Savage M (2015) Social Class in the 21st Century. Pelican Books, London Shukla N (ed) (2016) The good immigrant. Unbound, London\\

\noindent Steele C (2010) Whistling Vivaldi – and other clues to how stereotypes affect us. W. W. Norton \& Company, New York\\

\noindent Tatum B D (2003) Why are all the black kids sitting together in the cafeteria? Basic Books, New York\\

\noindent Todd, S (2021) Snakes and ladders: The Great British social mobility myth. Chatto \& Windus, London\\

\noindent Virdee S (2014) Racism, class and the racialized other. Bloomsbury Academic, London\\

\section*{In STEM Education}

Inequity and discrimination based on gender, race, ethnicity, (dis)ability, class, and sexuality observed and experienced in society is also documented to exist in educational contexts in both the United States and the United Kingdom. In 2015, the Council for Mathematical Sciences report on The Mathematical Sciences People Pipeline found that in the UK, female students make up about 40\% of undergraduate students in the Mathematical Sciences, but this number drops to 33\% in postgraduate students in the Mathematical Sciences. This drops again for Mathematics postgraduate students, of which 28\% identify as female. Black students make up 3\% of undergraduate students in the Mathematical Sciences, which is lower than the 6\% of the general undergraduate student population. The proportion of students with a reported disability is lower for postgraduate students (5\%) than undergraduate students in the Mathematical Sciences (8\%) and both of these are lower than in the general student population. Beyond university students, women made up only 15\% of all active authors in math, physics and computer science in the UK in 2010 (see Huang J., Gates A. J., Sinatra R., Barabási A.-L. (2020)). 

The US Equal Employment Opportunity Commission published its Annual Report: Women in STEM (FY2019) which states “Overall, women accounted for 29.3 percent of STEM federal workers. Science occupations had the most (49,546), while math occupations in the federal sector had the fewest number of women (6,469). There were significantly fewer women in Technology and Engineering than expected.”
A common metaphor to describe the phenomenon of more and more women dropping out of STEM as they advance through education and their career (at a higher rate than males), is that of a ‘leaky pipeline’.  The linked \href{https://www.catalyst.org/research/women-in-science-technology-engineering-and-mathematics-stem/}{Catalyst} website has extensive related data and information. 

It is generally acknowledged that the lack of diversity and inclusion stems less from overt exclusion happening at the level of student acceptance and recruitment, and more from deeper-seated societal pressures combined with implicit stereotypes and expectations. The longstanding \href{doi:10.1111/j.1949-8594.2002.tb18217.x}{Draw-a-Scientist} project both reveals and questions existing assumptions about the identity of a scientist as white, able, and male. One suggested explanation for this is that scientific epistemology represents an inherently white, masculine world view. Other existing research posits other possibilities, including: stereotypes of intellectual brilliance and implicit bias; absence of role models; greater accessibility of science curricula to men than women; lack of belonging, sexism, sexual harassment; and a notion that ``doing mathematics is doing masculinity.'' More generally, there is are questions of who can be a mathematical practitioner and what doing mathematics entails.\\ 

\noindent Bian L, Leslie S, Cimpian A (2018) Evidence of bias against girls and women in contexts that emphasize intellectual ability. American Psychologist 73:1139–1153. \url{https://doi.org/10.1037/amp0000427}\\

\noindent Bettinger E P, Long B T (2005) Do Faculty Serve as Role Models? The Impact of Instructor Gender on Female Students. The American Economic Review 95:152–157. \url{http://www.jstor.org/stable/4132808}\\

\noindent Blinkenstaff J C (2005) Women and science careers: Leaky pipeline or gender filter? Gender and Education 17:369–386. \url{https://doi.org/10.1080/09540250500145072}\\

\noindent Carlone H B, Johnson A (2007) Understanding the science experiences of successful women of color: Science identity as an analytic lens. JRST 44:1187-1218. \url{https://doi-org.ezproxy.st- andrews.ac.uk/10.1002/tea.20237}\\

\noindent Stephen J C, Shulamit K, Wendy M W (2023) Exploring Gender Bias in Six Key Domains of Academic Science: An Adversarial Collaboration. Psychological Science in the Public Interest,24:15-73. \url{https://doi.org/10.1177/15291006231163179}\\

\noindent Finson K D (2002) Drawing a Scientist: what we do and do not know after fifty years of drawings’, School Science and Mathematics 102:335–345. \url{doi:10.1111/j.1949- 8594.2002.tb18217.x}\\

\noindent Heffernan T (2022) Sexism, racism, prejudice, and bias: A literature review and synthesis of research surrounding student evaluations of courses and teaching. Assessment \& Evaluation in Higher Education 47:144–154. \url{https://doi.org/10.1080/02602938.2021.1888075}\\

\noindent Hill C, Corbett C, Rose R (2010) Why so few? Women in science, technology, engineering, and mathematics. American Association of University Women.\\
          
\noindent Huang J, Gates A J, Sinatra R, Barabási A (2020) Historical comparison of gender inequality in scientific careers across countries and disciplines. Proceedings of the National Academy of Sciences, USA, 117:4609–4616. \url{https://doi.org/10.1073/pnas.1914221117}\\

\noindent Faustino A C, Moura A Q, da Silva G H G, Muzinatti J L, Skovsmose O (2019) Microexclusion in Inclusive Mathematics Education. In: Kollosche D, Marcone R, Knigge M, Penteado M G, Skovsmose O (eds) Inclusive Mathematics Education. Springer, Cham. \url{https://doi.org/10.1007/978-3-030-11518-0_6}\\

\noindent Funk P, Iriberri N, Savio G (2019) Does scarcity of female instructors create demand for diversity among students? Evidence from an M-Turk experiment (CEPR Discussion Paper No. DP14190). SSRN. \url{https://ssrn.com/abstract=3504620}\\

\noindent Good C, Rattan A, Dweck C S (2012) Why do women opt out? Sense of belonging and women's representation in mathematics. Journal of Personality and Social Psychology 102:700–717. \url{https://doi.org/10.1037/a0026659}\\

\noindent Roy M, Guillopé C, Cesa M et al (2020) A Global Approach to the Gender Gap in Mathematical, Computing and Natural Sciences: How to Measure It, How to Reduce It?. International Mathematical Union. \url{https://zenodo.org/record/3882609}\\

\noindent Leslie S, Cimpian A, Meyer M, Freeland E (2015) Expectations of brilliance underlie gender distributions across academic disciplines. Science, 347:262–265. \url{https://doi.org10.1126/science.1261375}\\

\noindent Marcone R (2019) Who Can Learn Mathematics? in: Kollosche D, Marcone R, Knigge M, Penteado M G, Skovsmose O (eds) Inclusive Mathematics Education. Springer, Cham, pp. 41- 53\\

\noindent McGlaughlin S M, Knoop A J, Holliday G A (2005) Differentiating Students with Mathematics Difficulty in College: Mathematics Disabilities vs. No Diagnosis. Learning Disability Quarterly, 28:223–232. \url{https://doi.org/10.2307/1593660}\\

\noindent Mendick, H (2005) Mathematical stories: why do more boys than girls choose to study mathematics at AS-level in England? British Journal of Sociology of Education 26:235-251\\

\noindent Pezzoni M, Mairesse J, Stephan P, Lane J (2016) Gender and the publication output of graduate students: A case study. PLoS ONE 11:e0145146. \url{https://doi.org/10.1371/journal.pone.0145146}\\

\noindent Piatek-Jimenez K (2008) Images of Mathematicians: A New Perspective on the Shortage of Women in Mathematical Careers. ZDM – Mathematics Education 40:633–646. \url{https://doi- org.ezproxy.st-andrews.ac.uk/10.1007/s11858-008-0126-8}\\

\noindent Su R, Rounds J (2015) All STEM fields are not created equal: People and thing interests explain gender disparities across STEM fields. Frontiers in Psychology 6:189. \url{https://doi.org/10.3389/fpsyg.2015.00189}\\

\noindent Salinas P C, Bagni C (2017) Gender equality from a European perspective: Myth and reality. Neuron, 96:721–729. \url{https://doi.org/10.1016/j.neuron.2017.10.002}

\section*{Past and Present in Mathematics}
While sources above explore questions of equity, diversity, and inclusion in STEM educational contexts, this section focuses more narrowly on mathematics. References here vary widely in both intention and scope as they explore both historical development and contemporary factors that have contributed to persistent, systemic inequities and discrimination in the field. 

June Barrow-Green’s ‘The Historical Context of the Gender Gap in Mathematics,’ analyzes cultural attitudes to women involved in mathematics from the eighteenth century onwards. Similarly, Claire G. Jones’ \textit{Femininity, Mathematics and Science, 1880-1914} tracks the mechanisms of inclusions and exclusions in science and mathematics navigated by women at the turn of the twentieth century. In \textit{Masculinity and Science in Britain, 1831-1918}, Heather Ellis investigates ways in which mathematics and science gained a masculine identity. These three sources – and others on the list – highlight matters related to gender in mathematics. 
There are many additional factors that impact inequities and exclusions in the field. Texts by Joseph Dauben and Kapil Raj scrutinise approaches to the history of mathematics through an exclusively western lens. Others, including Theodore M. Porter and Ruth Schwartz Cowan, highlight how mathematical reasoning created and perpetuated individual categories of people. In another sense, the sources included here encourage the reader to consider the scope of mathematics, and to recognize its involvement particularly in the physical and invisible world around us.

Mathematical and technical tools have at times facilitated the construction of an exploitative reality. For example, Langdon Winner’s ‘Do Artifacts have Politics’ and subsequent book \textit{The Whale and the Reactor} explain how technological machines and constructions can contain and exercise power on society. More recent works invite an examination of everyday algorithms and how they might codify systemic inequities. Cathy O’Neil’s \textit{Weapons of Math Destruction} is at the forefront of this work, showcasing how mathematical algorithms perpetuate existing biases. Studies about the deployment of algorithms by police departments has been especially insightful; as seen in works by Virginia Eubanks and Sarah Brayne. Meanwhile, Joanna Radin’s “Digital Natives” shows how demographic data samples can be taken and reused in radically different contexts, and can perpetuate patterns of ‘settler colonialism’ in the process. As the influences of big data and algorithms grow, so, too, does the potential for them to perpetuate oppression. Above is a sample of the themes represented in this section of the reading list. We encourage the reader to dive in wherever their interest takes them.

\noindent Angwin J, Larson J, Mattu S, Kirchner L (2016) Machine bias: there’s software used across the country to predict future criminals. And it’s biased against blacks. ProPublica. \href{https://www.propublica.org/article/machine-bias-risk-assessments-in-criminal-sentencing}{https://www.propublica.org/article/machine-bias-risk-assessments-in-criminal-} \href{https://www.propublica.org/article/machine-bias-risk-assessments-in-criminal-sentencing}{sentencing}. Accessed 30 September 2023\\

\noindent Barrow-Green J (2019) The historical context of the gender gap in mathematics. In: Araujo C, Benkart G, Praeger C E, Tanbay B, (eds) World Women in Mathematics 2018: Proceedings of the First World Meeting for Women in Mathematics (WM)2. Springer, Cham, pp 129-144\\
 
\noindent Bouk D (2017) The history and political economy of personal data over the last two centuries in three acts. Osiris 32:85-106\\

\noindent Brayne S (2020) Coding inequality: how the use of big data reduces, obscures, and amplifies inequalities. In: Predict and surveil: data, discretion, and the future of policing. Oxford University Press, Oxford, pp 100-117\\

\noindent Carter R G S (2006) Of things said and unsaid: power, archival silences, and power in silence. Archivaria 61:215-233.\\

\noindent Cheney-Lippold J (2017) We are data: algorithms and the making of our digital selves. NYU Press, New York\\

\noindent Cowan R S (1972) Francis Galton’s statistical ideas: the influence of eugenics. Isis 63:509-528. 10.1086/351000\\

\noindent Dauben J W (2021) Anachronism and incommensurability: words, concepts, contexts, and intentions. In: Guicciardini N (ed) Anachronism in the history of mathematics: essays on the historical interpretation of mathematical texts. Cambridge University Press, Cambridge, pp 307-357\\

\noindent Ellis H (2017) Masculinity and science in Britain. Palgrave Macmillan, London\\

\noindent Eubanks V (2019) Automating inequality: how high-tech tools profile, police and punish the poor. St Martin’s Press, New York\\

\noindent Jones C G (2009) Femininity, mathematics and science, 1880-1914. Palgrave Macmillan, London\\

\noindent Jones C G, Martin A E, Wolf A (eds) (2022) The palgrave handbook of women in science since 1660. Palgrave Macmillan, Cham
Michel N (2021) Mathematical selves and the shaping of mathematical modernism: conflicting epistemic ideals in the emergence of enumerative geometry (1864–1893). Isis 112:68-92. \url{https://doi.org/10.1086/713831}\\

\noindent Milam E L, Nye R A (2015) An introduction to scientific masculinities. Osiris 30:1-14.
\url{https://doi.org/10.1086/682953}\\
  
\noindent Muhammad K G (2010) Saving the nation: the racial data revolution and the negro problem. In: The condemnation of blackness: race, crime, and the making of modern urban America. Harvard University Press, Cambridge, pp 15-34\\

\noindent Murray M (2001) Women becoming mathematicians: creating a professional identity in post- World War II America. The MIT Press, Cambridge\\

\noindent Niskanen K, Barany M J (2021) Gender, embodiment and the history of the scholarly persona: incarnations and contestations. Palgrave Macmillan, Cham\\

\noindent Noble S (2018) Algorithms of Oppression: How Search Engines Reinforce Racism. NYU Press\\

\noindent O’Neil C (2017) Weapons of math destruction: how big data increases inequality and threatens democracy. Crown Books, London\\

\noindent Porter T M (1995) The political philosophy of quantification. In: Trust in numbers: the pursuit of objectivity in science and public life. Princeton University Press, Princeton, pp 72-86\\

\noindent Porter T M (2018) Genetics in the madhouse: the unknown history of human heredity. Princeton University Press, Princeton\\

\noindent Powles J, Nissenbaum H (2018) The seductive diversion of ‘solving’ bias in artificial intelligence. OneZero. \url{https://onezero.medium.com/the-seductive-diversion-of-solving- bias-in}\\\url{-artificial-intelligence-890df5e5ef53}. Accessed 30 September 2023\\

\noindent Radin J (2017) ‘Digital natives’: how medical and indigenous histories matter for big data. Osiris 32:43-64. \url{https://doi.org/10.1086/693853}\\

\noindent Raj K (2007) Relocating modern science: circulation and the construction of knowledge in South Asia and Europe 1650-1900. Palgrave Macmillan, New York\\

\noindent Raj K (2013) Beyond postcolonialism ... and positivism: circulation and the global history of science. Isis 104:337-347. \url{https://doi.org/10.1086/670951}\\

\noindent Rittenberg C J, Tanswell F S, Van Bendegem J P (2018) Epistemic injustice in mathematics. Synthese 197:3875-3904. \url{https://doi.org/10.1007/s11229-018-01981-1}\\

\noindent Robson E, Stedall J (eds) (2009) The Oxford handbook of the history of mathematics. Oxford University Press, Oxford\\
     
\noindent Sasada M, Bannaj K (2021) The situation of gender equality in mathematics in Japan. \url{https://www.math.keio.ac.jp/~bannai/Report_MathGender_en.pdf.} Accessed 30 September 2023\\

\noindent Schappacher N (2022) Framing global mathematics: the international mathematical union between theorems and politics. Springer, Cham\\

\noindent Seltzer W, Anderson M (2001) The dark side of numbers: the role of population data systems in human rights abuses. Social Research 68:481-513\\

\noindent Tanswell F S, Kidd I J (2021) Mathematical practice and epistemic virtue and vice. Synthese 199:407-426. \url{https://doi.org/10.1007/s11229-020-02664-6}

\noindent Winner L (1980) Do artifacts have politics?. Daedalus 109:121-136\\

\noindent Winner L (2020) The whale and the reactor: a search for limits in an age of high technology. University of Chicago Press, Chicago\\

\section*{Alternative Visions and Approaches}

Given the above articulations of various ways mathematical education faces challenges in this social context, this section contains materials with suggestions on how one might approach mathematics and mathematics education differently. Central to such efforts is the realization that mathematics is a social practice, which has been shaped by (and reciprocally also has shaped) historical and social trends: it is not acontextual or independent from human experience. Rather than something existing independently from humans, mathematical practice is continuously constructed and re-constructed by humans. Glas (2006) investigates how mathematics can be both objective and a historically contingent practice.
 
The acknowledgment of mathematics as dependent on human context is essential within a discussion of DEI/EDI because it counters the pervasive notion that mathematics is  based on intuition, possessed by some and not others. This myth serves to justify arguments which posit only certain people as capable of doing/understanding mathematics. Further, once mathematics is considered as a social practice, it becomes possible to acknowledge that different societies have different ways of doing mathematics, without positing a hierarchy between them. This insight was central to the development of ethnomathematics, a field pioneered by Ubiratan D’Ambrosio, who critiqued the separation of mathematics from questions of social justice and cultural issues. Many of the articles included in this section thus argue that mathematics education should address questions of justice and ethics, and aim to include students’ backgrounds in the learning of mathematics. 
 
From such a point-of-view, the history of mathematics is not only essential to the realization that mathematics has been constantly changing through societal interaction and  different ways of doing mathematics have always existed, but it also becomes a central feature of mathematics education. While some argue that the inclusion of history in mathematics education serves to deepen the knowledge/understanding of students and offers different ways of learning and approaching mathematics (for such approaches see the section on using primary sources for mathematics education), others see it as a valuable end in itself. Indeed, a study of the history of mathematics confronts students with different ways of doing mathematics. This can both interrogate existing assumptions about mathematical knowledge production and invite students to position themselves not as passive recipients, but as meaning-makers in the classroom who are navigating their identity in relation to mathematics.   

There are a number of websites designed to increase visibility of diverse role models in mathematics. The website that is now
\href{https://www.mathad.com/home}{Mathematicians of the African Diaspora} was started in 1997 by Scott Williams to promote the contributions to mathematical research from members of the African diaspora. The name of the \href{https://mathematicallygiftedandblack.com/}{Mathematically Gifted \& Black} website comes from a song sung by Nina Simone and co-written by Weldon Irvine. It was created in 2016 and features the accomplishments of Black scholars in the mathematical sciences. Marie A. Vitulli created The \href{https://pages.uoregon.edu/wmnmath/biographies.html}{Women in Math Project} in 1997 to make accessible some resources and information for and about women in mathematics. It is still hosted by the University of Oregon. As part of an ongoing project hosted by Agnes Scott College, \href{https://mathwomen.agnesscott.org/women/women.htm}{Biographies of Women mathematicians} spotlights achievements in mathematics from both contemporary and historical women.  The \href{https://www.sacnas.org/sacnas-biography-project}{SACNAS} biography project has an online archive of stories by and about Chicano/Hispanic and Native American scientists with advanced degrees. Spectra, The Association for LBGTQ+ Mathematicians hosts \href{http://lgbtmath.org/People.html}{Out Lists} to provide resources and support for interested mathematicians.

Additional sources suggest other ways that careful and thoughtful implementation of research in the history of mathematics can contribute to the discipline developing in a direction of greater inclusivity. Work in the history of mathematics can provide a platform from which to view the discipline critically and many of these works aim to alter our view of mathematical knowledge, its production and people. They have expanded the discussion to include previously ignored contributions, to take a global view of mathematical practice, and to challenge existing notions about the mathematical professions. Such resources can assist in highlighting past incidents of oppression and discrimination, while uncovering some roots and causes of persistent legacies today. Research in and resources from the history of mathematics can also be effective tools to help build a more inclusive discipline. 

Additionally, there exists a body of work about using primary historical sources to teach technical mathematical content. This can facilitate differently accessible learning and also create connections beyond the standard curricular content.\\ 

\noindent Anderton L, Wright D (2012) We could all be having so much more fun!: A case for the history of mathematics in education. Journal of Humanistic Mathematics 2:88-103. 10.5642/jhummath.201201.08\\

\noindent Azzouni J (2006) How and why mathematics is unique as social practice. In: Hersh R (ed) 18 Unconventional essays on the nature of mathematics. Springer, New York, pp 201-219\\
     
\noindent Chemla K (2018) How has one, and how could have one, approached the diversity of mathematical cultures? In: Mehrmann V, Skutella M (eds) European Congress of Mathematics: Berlin, July 18-22, 2016. European Mathematical Society Press, Berlin, pp 1-61\\

\noindent Cliffe E, Rowlett P (eds) (2012) Good practice in inclusive curricula in the mathematical sciences. Maths, Stats \& OR Network. \url{https://www.mathcentre.ac.uk/resources/uploaded/inclusivecurricula.pdf}. Accessed 30 September 2023\\

\noindent D’Ambrosio U (2001) Mathematics and peace: a reflection on the basis of Western civilization. Leonardo 34: 327-332\\

\noindent D’Ambrosio U (2016) An overview of the history of ethnomathematics. In: Rosa M, D’Ambrosio U, Orey D C, Shirley L, Alangui W V, Palhares P, Gavarrete M E (eds) Current and future perspectives of ethnomathematics as a program. Springer, Cham, pp 5-10\\

\noindent Fried M (2008) History of mathematics in mathematics education: a Saussurean perspective. The Mathematics Enthuisiast 5:185-198. \url{https://doi.org/10.54870/1551-3440.1100}\\

\noindent Glas E (2006) Mathematics as objective knowledge and as human practice. In: Hersh R (ed) 18 Unconventional essays on the nature of mathematics. Springer, New York, pp 289-303\\

\noindent Green B (2021) Data science as political action: grounding data science in a politics of justice. Journal of Social Computing 2:249-265. \url{https://doi.org/10.48550/arXiv.1811.03435}\\

\noindent Izmirli I H (2011) Pedagogy on the ethnomathematics-epistemology nexus: a manifesto. Journal of Humanistic Mathematics, 1:27-50. 10.5642/jhummath.201102.04\\

\noindent Jankvist U T (2009) A categorization of the ‘whys’ and ‘hows’ of using history in mathematics education. Educational Studies in Mathematics 71:235-261. \url{https://doi.org/10.1007/s10649- 008-9174-9}\\

\noindent Mayes-Tang S (2019) Telling women’s stories: a resource for college mathematics instructors. Journal of Humanistic Mathematics 9:78-92. 10.5642/jhummath.201902.07\\

\noindent Nunez R (2006) Do real numbers really move? Language, thought and gesture: the embodies cognitive foundations of mathematics. In: Hersh R. (ed) 18 Unconventional essays on the nature of mathematics. Springer, New York, pp 160-181\\

\noindent Peck F (2021) Towards anti-deficit education in undergraduate mathematics education: how deficit perspectives work to structure inequality and what can be done about it. PRIMUS 31:940-961. \url{https://doi.org/10.1080/10511970.2020.1781721}\\

\noindent Radford L, Santi G (2022) Learning as a critical encounter with the other: prospective teachers conversing with the history of mathematics. ZDM Mathematics Education 54:1479-1492. \url{https://doi.org/10.1007/s11858-022-01393-z}\\

\noindent Schulman B (2002) Is there enough poison gas to kill the city?: the teaching of ethics in mathematics classes. The College Mathematics Journal 33:118-125. \url{https://doi.org/10.2307/1558994}\\

\noindent Simic-Muller K (2023) Noticing and wondering to rehumanize mathematics classrooms. PRIMUS 33:431-443. \url{https://doi.org/10.1080/10511970.2022.2073624}\\

\noindent Skovomose O (2019) Inclusions, meetings and landscapes. In: Kollosche D, Marcone R, Knigge M, Penteado M G, Skovomose O (eds) Inclusive mathematics education: state-of-the-art research from Brazil and Germany. Springer, Cham, pp 71-84\\

\noindent Spindler R (2022) Foundational mathematical beliefs and ethics in mathematical practice and education. Journal of Humanistic Mathematics, 12:49-71. 10.5642/jhummath.GOSN2205\\

\noindent Tanswell F S, Rittberg C J (2020) Epistemic injustice in mathematics education. ZDM 52:1199- 1210. \url{https://doi.org/10.1007/s11858-020-01174-6}\\

\noindent Yeh C, Rubel L (2020) Queering mathematics: disrupting binary oppositions in mathematics pre-service education. In: Radakovic N, Jao L (eds) Borders in mathematics pre-service teacher education. Springer, Cham, pp 227-243\\

\noindent Webb S D (2011) Accessibility of university mathematics. MSOR Connections 11:42-46\\

\subsection*{Using Primary Sources to teach Mathematics}

\noindent Barnett J H (2014) Learning Mathematics via Primary Historical Sources: Straight From the Source’s Mouth, PRIMUS, 24:722-736, \url{https://doi.org/10.1080/10511970.2014.899532}\\

\noindent Barnett J H, Lodder J, Pengelley D (2016) Teaching and Learning Mathematics From Primary Historical Sources, PRIMUS, 26:1-18, \url{https://doi.org/10.1080/10511970.2015.1054010}\\

\noindent Chemla K (2012) Using documents from Ancient China to Teach Mathematical Proof. In Hanna G, de Villiers M (eds) Proof and Proving in Mathematics Education. Springer, Dordrecht, pp. 423-429\\

\noindent Clark K, Kjeldsen T H, Schorcht S, Tzanakis C (eds) (2018) Mathematics, Education and History: Towards a Harmonious Partnership. Springer, Cham.\\
     
\noindent Clark, K M, Can C, Barnett J H et al (2022) Tales of research initiatives on university-level mathematics and primary historical sources. ZDM – Mathematics Education 54:1507–1520. \url{https://doi.org/10.1007/s11858-022-01382-2}\\

\noindent Fried M N (2001) Can Mathematics Education and History of Mathematics Coexist? Science \& Education 10:391-408. \url{https://doi.org/10.1023/A:1011205014608}\\

\noindent Fried M N (2003) Humanistic Mathematics as Mathematics for All. Humanistic Mathematics Network Journal 27:15. \url{https://doi.org/10.5642/hmnj.200401.27.15}\\

\noindent Fried M N (2002) Edmond Halley and Apollonius: Second-Order Historical Knowledge in Mathematics Education. ZDM – Mathematics Education 54:1435-1447. \url{https://doi.org/10.1007/s11858-022-01391-1}\\

\noindent Lodder J (2014) Networks and Spanning Trees: The Juxtaposition of Prüfer and Borůvka. PRIMUS, 24:737-752. \url{https://doi.org/10.1080/10511970.2014.896835}\\

\noindent Mayfield B (2014) Weaving History Through the Major. PRIMUS 24:669-683. \url{https://doi.org/10.1080/10511970.2014.900158}\\

\section*{Conclusion }

One purpose of this source list is to facilitate an aim articulated in the recent QAA Benchmark Statement for MSOR: “values of EDI should permeate the curriculum and every aspect of the learning experience to ensure the diverse nature of society in all its forms is evident.”  Although the collection of sources listed here is only a sample of available resources, the list has been assembled as an aid to make this work less daunting. It is not possible to provide a single document that will target every particular concern or priority given the wide range of contexts in which mathematical practitioners work. Specific implementation in a particular context will necessarily involve personal effort and perhaps attempts and revisions. This collection of sources nonetheless aims to give a broad range of documents that can support DEI/EDI efforts. 
Particularly important to this resource is the inclusion of work from the history of mathematics alongside texts that document DEI/EDI issues in broad cultural contexts as well as more specific educational contexts. The history of mathematics is especially well positioned to facilitate a move within the discipline of mathematics towards greater inclusivity. Historical studies can not only illuminate roots and causes of oppressive practices, but also challenge dominant narratives of who mathematics includes. The QAA statement states “it is highly desirable that students encounter a wide range of role models within higher education.” While the community works to elevate an array of identifiable role models, sources listed here can perhaps help ameliorate challenging circumstances. In these ways and others, historical scholarship provides resources to help set a brighter course for the future of the mathematical sciences. 

\end{document}